\documentclass[11pt]{amsart}

\usepackage{amscd}
\usepackage{amsmath}
\usepackage{graphicx}
\usepackage{amsfonts}
\usepackage{amssymb}
\begin{document}
\title[Maps preserving peripheral spectrum]
{Maps preserving peripheral spectrum of generalized products
of operators}

\author{Wen zhang}
\address[Wen Zhang ]{Department of
 Mathematics, Shanxi University, Taiyuan, Shanxi,
030006, P. R. China}\email{wenzhang1314gw@163.com}

\author{Jinchuan Hou}
\address[Jinchuan Hou]{Department of Mathematics, Taiyuan
University of Technology, Taiyuan 030024, P. R. of China; Department
of
 Mathematics, Shanxi University, Taiyuan, Shanxi,
030006, P. R. China} \email{houjinchuan@tyut.edu.cn;
jinchuanhou@aliyun.com.cn}

\thanks{This work is partially  supported by National Natural Science Foundation
of China (No.11171249,   11271217).}

\begin{abstract}
Let $\mathcal{A}_1$ and $\mathcal{A}_2$ be standard operator
algebras on complex Banach spaces $X_1$ and $X_2$, respectively. For
$k\geq2$, let $(i_1,...,i_m)$ be a sequence with terms chosen from
$\{1,\ldots,k\}$, and assume that at least one of the terms in
$(i_1,\ldots,i_m)$ appears exactly once. Define the generalized
product $T_1* T_2*\cdots* T_k=T_{i_1}T_{i_2}\cdots T_{i_m}$ on
elements in $\mathcal{A}_i$.  Let
$\Phi:\mathcal{A}_1\rightarrow\mathcal{A}_2$ be a map with the range
  containing all operators of rank at most two. We show  that $\Phi$ satisfies that
$\sigma_\pi(\Phi(A_1)*\cdots*\Phi(A_k))=\sigma_\pi(A_1*\cdots* A_k)$
for all $A_1,\ldots, A_k$, where $\sigma_\pi(A)$ stands for the
peripheral spectrum of $A$, if and only if $\Phi$ is an  isomorphism
or an anti-isomorphism multiplied by an $m$th root of unity, and the
latter case occurs only if the generalized product is quasi-semi
Jordan. If $X_1=H$ and $X_2=K$ are complex Hilbert spaces, we
characterize also maps preserving the peripheral spectrum of the
skew generalized products, and prove that such maps are of the form
$A\mapsto cUAU^*$ or $A\mapsto cUA^tU^*$, where
$U\in\mathcal{B}(H,K)$ is a unitary operator, $c\in\{1,-1\}$.

\end{abstract}
\maketitle

\section{Introduction}

Linear maps between Banach algebras which preserve the spectrum are
extensively studied in connection with a longstanding open problem
due to Kaplansky on invertibility preserving linear maps (\cite{A1,
AM, CH1,CH2,CH3,HC, JS, K} and the references therein). Recently,
the study of spectrum preservers without linearity or additivity
assumption  also attracted attentions of researchers. One of
interesting topics of this kind concerns the spectrum of products.
In \cite{M1}, Moln$\acute{a}$r characterized surjective maps $\Phi$
on bounded linear operators acting on a Hilbert space preserving the
spectrum of the product of operators, i.e., $AB$ and
$\Phi(A)\Phi(B)$ always have the same spectrum. This similar
question was studied by Huang and Hou in \cite{HH} by replacing the
spectrum by several spectrual functions such as the left spectrum,
spectral boundary, etc.. Hou, Li and Wong \cite{HLW1} studied
further the maps $\Phi$ between certain operator algebras preserving
the spectrum of a generalized product $T_1* T_2*\cdots* T_k$ of low
rank operators. Namely, for all operators $T_1, T_2,\cdots,T_k$ of
low rank the spectra of $T_1* T_2*\cdots* T_k$ and of
$\Phi(T_1)*\Phi(T_2)*\cdots* \Phi(T_k)$ are equal. The generalized
product   is defined as following.

{\bf Definition 1.1.} {\it Fix a positive integer $k$ and a finite
sequence $(i_1,i_2,\ldots,i_m)$ such that
$\{i_1,i_2,\ldots,i_m\}=\{1,2,\ldots,k\}$ and there is an $i_p$ not
equal to $i_q$ for all other $q$. A generalized product for
operators $T_1,\ldots,T_k$ is defined by   $$T_1* T_2*\cdots*
T_k=T_{i_1}T_{i_2}\cdots T_{i_m}, \eqno(1.1)$$ $m$ is called the
width of the generalized product.

Furthermore, if $(i_1,i_2,\ldots,i_m)$ is symmetrical with respect
to $i_p$, we say that $T_1* T_2*\cdots* T_k$ is a   generalized semi
Jordan product; if $$(i_{p+1}, i_{p+2},\ldots,i_{m-1},i_m, i_1,
i_2,\ldots, i_{p-1})=(i_{p-1},i_{p-2},\ldots, i_2,i_1,
i_m,i_{m-1},\ldots, i_{p+2},i_{p+1}), \eqno(1.2)$$ we say that $T_1*
T_2*\cdots* T_k$ is a generalized quasi-semi Jordan product.}

 Evidently, this
definition of generalized product covers the usual product $T_1T_2$,
Jordan semi-triple $BAB$ and the triple one:
$\{T_1,T_2,T_3\}=T_1T_2T_3$, etc.; the definition of generalized
semi Jordan product cover the Jordan semi-triple $BAB$ and the
product like $T_1*T_2*T_3= T_2T_3^2T_1T_3^2T_2$; the definition of
generalized quasi-semi Jordan product covers the   products like
$A*B=B^rAB^s$ and $ T_1*T_2*T_3=T_2T_3^2T_1T_3^2T_2^2T_3^3T_2$.

 Let
$\mathcal{B}(X)$ be the Banach algebra of all bounded linear
operators on a complex Banach space $X$.  Recall that a standard
operator algebra ${\mathcal A}$ on a complex Banach space $X$
usually stands for  a closed subalgebra of $\mathcal{B}(X)$
containing the ideal of all finite rank operators  and the identity
$I$ on $X$.  However, in the present paper, we do not assume that
${\mathcal A}$ contains the identity operator, or that it is closed.

Denote by $\sigma(T)$ and $r(T)$ the spectrum and the spectral
radius of $T\in \mathcal{B}(X)$, respectively. The peripheral
spectrum of $T$ is defined by
$$\sigma_\pi(T)=\{z\in\sigma (T):|z|=r(T)\}.$$ Since $\sigma(T)$ is
compact, $\sigma_\pi(T)$ is a well-defined non-empty set and is an
important spectral function. Also observe that
$\sigma_\pi(TS)=\sigma_\pi(ST)$ holds for any $T,S\in{\mathcal
B}(X)$.

In \cite{TL}, Tonev and Luttman studied maps preserving peripheral
spectrum of the usual operator products on standard operator
algebras. It was proved that, if such a map is surjective, then it
must be a positive or negative multiple of an isomorphism or an
anti-isomorphism. They studied also the corresponding problems in
uniform algebras  (see \cite{LLT,LT}). Recently, Takeshi and Dai
\cite{TD} generalized the result in \cite{LT}, and characterized
surjective maps $\phi$ and $\psi$ satisfying
$\sigma_\pi(\phi(T)\psi(S))=\sigma_\pi(TS)$ on standard operator
algebras. The surjective maps between standard operator algebras on
Hilbert spaces that preserve the peripheral spectrum of skew
products $T^*S$ of operators was also characterized in \cite{TD}.
Cui and Li studied in \cite{CL} the maps preserving peripheral
spectrum of Jordan products of operators $AB+BA$ on standard
operator algebras. They show that, if the range of such a map
contains all operators of rank at most 2, then it is an isomorphism
or an anti-isomorphisms multiplied by $1$ or $-1$. A
characterization of maps preserving the peripheral spectrum of skew
Jordan products $AB^*+B^*A$ was also given in \cite{CL}. In
\cite{ZH1} the  maps preserving peripheral spectrum of Jordan
semi-triple products $BAB$ of operators is characterized.

Motivated by the above results, we consider the question of
characterizing the maps preserving the peripheral spectrum of the
generalized products of operators defined in Eq.(1.1). In fact, the
purpose of this paper is manifold. Firstly, we characterize maps
between standard operator algebra on Banach spaces preserving
peripheral spectrum of the generalized product of operators under a
mild assumption that  the range of the map contains all operators of
rank at most two. Let $\mathcal{A}_1$ and $\mathcal{A}_2$ be
standard operator algebras on complex Banach spaces $X_1$ and $X_2$,
respectively. Assume that
$\Phi:\mathcal{A}_1\rightarrow\mathcal{A}_2$ is a map the range of
which contains all operators of rank at most two. We show that
$\sigma_\pi(\Phi(A_1)*\cdots*\Phi(A_k))=\sigma_\pi(A_1*\cdots* A_k)
$ holds for all $A_1,A_2,\ldots, A_k\in{\mathcal A}_1$ if and only
if either there exist a scalar $\lambda\in\mathbb{C}$ with
$\lambda^m=1$ and an invertible operator $T\in\mathcal{B}(X_1,X_2)$
such that $\Phi(A)=\lambda TAT^{-1}$ for all $A\in\mathcal{A}_1$; or
there exists a scalar $\lambda\in\mathbb{C}$ with $\lambda^m=1$ and
an invertible operator $T\in\mathcal{B}(X^*_1,X_2)$ such that
$\Phi(A)=\lambda TA^*T^{-1}$ for all $A\in\mathcal{A}_1$. In the
last case, the spaces $X_1$ and $X_2$ must be reflexive,
$A_1*\cdots* A_k$   a general quasi-semi Jordan product or $k=2$
(see Theorem 2.1). Particularly, if the generalized product is not
semi Jordan and $k\geq 3$, then $\Phi$ preserves the peripheral
spectrum of the generalized product   if and only if $\Phi$ is an
isomorphism multiplied by an $m$th root of 1.
 Secondly, we characterize maps
preserving the peripheral spectrum of the skew generalized product
of operators on Hilbert space. As expected such maps are
$*$-isomorphism or $*$ anti-isomorphism; or, in the case $m$ is
even, $*$-isomorphism or $*$ anti-isomorphism multiplied by   $-1$
(see Theorem 3.1).

Throughout this paper,  $X$ stands for  complex Banach spaces of any
dimension. Denote by $X^*$ the dual space of $X$ and by
$\mathcal{B}(X)$ the Banach algebra of all bounded linear operators
on $X$. For $A\in\mathcal{B}(X)$, $A^*$ denotes the adjoint operator
of $A$. For nonzero $x\in X $ and $f\in X^*$,  $x\otimes f$ is the
rank one operator defined by $(x\otimes f)z = f(z)x$ for every $z\in
X$. We often use $\langle x,f\rangle$ for $f(x)$, the value of $f$
at $x$. For $A\in\mathcal{B}(X)$, $\ker(A)$ and ${\rm ran}(A)$
denote respectively the kernel and the range of $A$, while rank$(A)$
stands for the rank of $A$, that is, the dimension of ran$(A)$. Let
$\mathbb{C}$ and $\mathbb{R}$ denote respectively the complex field
and real field as usual.

\section{Generalized products of operators on Banach space}

In this section, we  study maps between standard operator algebras
on complex Banach spaces preserving peripheral spectrum of the
generalized products of operators. It is clear that every
isomorphism between standard operator algebras preserves the
peripheral spectrum of any generalized product of operators. Recall
that a Jordan isomorphism $\Phi:
\mathcal{A}_1\rightarrow\mathcal{A}_2$ is either a spacial
isomorphism or a spacial anti-isomorphism. In this case, for any
generalized semi Jordan product  $A_1*\cdots* A_k$,
$\sigma_\pi(\Phi(A_1)* \cdots*\Phi(A_k))=\sigma_\pi(A_1*\cdots*
A_k)$ holds for all $A_1, \ldots,A_k$.  Our main result below gives
a characterization of the maps between standard operator algebras
that preserve the peripheral spectrum of generalized products of
operators.

{\bf Theorem 2.1.} {\it Consider the product $T_1*\cdots* T_k$
defined in Definition 1.1 with width $m$. Assume that
$\Phi:\mathcal{A}_1\rightarrow\mathcal{A}_2$ is a map the range of
which contains all operators of rank at most two. Then $\Phi$
satisfies
$$\sigma_\pi(\Phi(A_1)*\cdots*\Phi(A_k))=\sigma_\pi(A_1*\cdots*
A_k) \eqno(2.1)$$   for all $A_1,A_2,\ldots, A_k\in{\mathcal A}_1$
if and only if   one of the following conditions holds. }

(1) {\it There exist a scalar $\lambda\in\mathbb{C}$ with $\lambda^m=1$
and an invertible operator $T\in\mathcal{B}(X_1,X_2)$ such that
$\Phi(A)=\lambda TAT^{-1}$ for all $A\in\mathcal{A}_1$.}

(2) {\it  There exists a scalar $\lambda\in\mathbb{C}$ with
$\lambda^m=1$ and an invertible operator
$T\in\mathcal{B}(X^*_1,X_2)$ such that $\Phi(A)=\lambda TA^*T^{-1}$
for all $A\in\mathcal{A}_1$. In this case, the spaces $X_1$ and
$X_2$ are reflexive, $A_1*\cdots* A_k$ is a generalized quasi-semi
Jordan product.}

By Theorem 2.1, if $X_i$ is not reflexive for some $i=1,2$, or, if
the generalized product is not quasi-semi Jordan, then $\Phi$
satisfies Eq.(2.1) if and only if $\Phi$ is an isomorphism
multiplied by an $m$th root of 1.

 To prove  Theorem 2.1, we first consider the special
case of $k=2$. Thus there exists nonnegative integers $r,s$ with
$r+s=m-1\geq 1$ such that $A_1*A_2=A_2^rA_1A_2^s$.

{\bf Theorem 2.2.} {\it Let $\mathcal{A}_1$ and $\mathcal{A}_2$ be
standard operator algebras on complex Banach spaces $X_1$ and $X_2$,
respectively. Assume that
$\Phi:\mathcal{A}_1\rightarrow\mathcal{A}_2$ is a map the range of
which contains all operators of rank at most two, and $r,s$ are
nonnegative integers with $r+s\geq 1$. Then $\Phi$ satisfies
$$\sigma_\pi(B^rAB^s)=\sigma_\pi(\Phi(B)^r\Phi(A)\Phi(B)^s)\quad\mbox{\it for all}\
 A,B\in\mathcal{A}_1 \eqno(2.2)$$ if and only if  one of the following two
 statements holds.}

 (1) {\it There exist a scalar $\lambda\in\mathbb{C}$ with $\lambda^m=1$
and an invertible operator $T\in\mathcal{B}(X_1,X_2)$ such that
$\Phi(A)=\lambda TAT^{-1}$ for all $A\in\mathcal{A}_1$.}

(2) {\it There exists a scalar $\lambda\in\mathbb{C}$ with
$\lambda^m=1$ and an invertible operator
$T\in\mathcal{B}(X^*_1,X_2)$ such that $\Phi(A)=\lambda TA^*T^{-1}$
for all $A\in\mathcal{A}_1$. In this case, the spaces $X_1, X_2$ are
reflexive.}

It is interesting to remark here that above results for the
peripheral spectrum are some what different from  the corresponding
results for the spectrum. In fac, $\Phi$ satisfies $\sigma
(\Phi(A_1)*\cdots*\Phi(A_k))=\sigma (A_1*\cdots* A_k)$ for all
$A_1,A_2,\ldots, A_k\in{\mathcal A}_1$ (resp.
$\sigma(B^rAB^s)=\sigma(\Phi(B)^r\Phi(A)\Phi(B)^s)$ for all
$A,B\in{\mathcal A}_1$)   if and only if either $\Phi$ has the form
(1) of Theorem 2.1, or

(2$^\prime$)  There exists a scalar $\lambda\in\mathbb{C}$ with
$\lambda^m=1$ and an invertible operator
$T\in\mathcal{B}(X^*_1,X_2)$ such that $\Phi(A)=\lambda TA^*T^{-1}$
for all $A\in\mathcal{A}_1$. In this case, the spaces $X_1$ and
$X_2$ are reflexive; moreover, $A_1*\cdots* A_k$ is a generalized
semi Jordan product (resp. $r=s$) whenever there exist left (or
right) invertible element in ${\mathcal A}_1$ that is not
invertible.

Now we apply Theorem 2.2 to prove Theorem 2.1.

{\bf Proof of Theorem 2.1.} We may assume that $\dim X_1\geq 2$. For
the ``if" part, (1)$\Rightarrow$ Eq.(2.1) is obvious;
(2)$\Rightarrow$ Eq.(2.1) because $\sigma_\pi(AB)=\sigma_\pi(BA)$
and the generalized product is quasi-semi Jordan. To check the
``only if" part, assume that $\Phi$ satisfies Eq.(2.1). Consider the
special case of generalized product $A_1*A_2*\cdots
*A_k$  with $A_{i_p}=A$ and all other $A_{i_q}=B$, one sees that
$\Phi$ satisfies Eq.(2.2). By Theorem 2.2, $\Phi$ has the form (1)
or the form (2) in Theorem 2.2.

To complete the proof, by Theorem 2.2, we need only to show that
$\Phi$ takes the form (2) will imply that  $A_1*A_2*\cdots
*A_k$ is a generalized quasi-semi Jordan product. Note that we always have $\sigma_\pi(AB)=\sigma_\pi(BA)$.
Then, as $\Phi(A)=\lambda TA^*T^{-1}$ for any $A\in{\mathcal A}_1$,
we see that
$$\begin{array}{rl} \sigma_\pi(A_{i_1}A_{i_2}\cdots A_{i_m})=& \sigma_\pi(A_1*\cdots*
A_k)=\sigma_\pi(\Phi(A_1)*\cdots*\Phi(A_k))\\
= & \lambda^m \sigma_\pi(A_{i_1}^*A_{i_2}^*\cdots A_{i_m}^*)
=\sigma_\pi(A_{i_m}A_{i_{m-1}}\cdots A_{i_1})\end{array}$$ holds for
all $A_1,A_2,\ldots, A_k\in{\mathcal A}_1$.  Thus, with  $i_p$ the
same as in Definition 1.1, one has
$$\begin{array}{rl} &\sigma_\pi(A_{i_p}A_{i_{p+1}}\cdots A_{i_m}A_{i_1}A_{i_2}\cdots A_{i_{p-1}})=\sigma_\pi(A_{i_1}A_{i_2}\cdots A_{i_m})\\
= & \sigma_\pi(A_{i_m}A_{i_{m-1}}\cdots
A_{i_1})=\sigma_\pi(A_{i_p}A_{i_{p-1}}\cdots A_{i_1}A_{i_m}\cdots
A_{i{p+1}})\end{array}$$ holds for all $A_1,A_2,\ldots,
A_k\in{\mathcal A}_1$. This implies, by a similar argument in
\cite[Theorem 3.2]{HLW1}, that Eq.(1.2) holds and hence, the
generalized product is quasi-semi Jordan. \hfill$\Box$

To prove Theorem 2.2, as one may expect, we will show that $\Phi$ is
linear and  preserves rank one operators in both directions.  The
following lemma is crucial, which gives a characterization of rank
one operators in terms of the peripheral spectrum of the generalized
products.

{\bf Lemma 2.3.} {\it Suppose $r$ and $s$ are nonnegative integers
such that $r+s\geq 1$. Let $A$ be a nonzero operator on a complex
Banach space $X$ of dimension at least two. Then the following
conditions are equivalent.}

(1) {\it $A$ is of rank one.}

(2) {\it For any $B\in\mathcal{A}$, $ \sigma_\pi(B^rAB^s)$ is a singleton.}

(3) {\it For any $B\in\mathcal{A}$ with {\rm rank}$(B)\leq2$,
$\sigma_\pi(B^rAB^s)$ is a singleton.}

{\bf Proof}. The implications $(1)\Rightarrow(2)\Rightarrow(3)$ are
clear.

To prove $(3)\Rightarrow(1)$, we consider the contrapositive. Since
the case $r+s=1$ is easily checked, we assume in the rest of the
proof that $r+s\geq 2$.

Suppose (3) holds but (1) is not true, i.e., $A$ has rank at least
two. Then there exist linearly independent vectors $x_1,x_2\in X$
such that $\{Ax_1,Ax_2\}$ is a linearly independent set. Fix such
$x_1$ and $x_2$. We complete the proof by considering the following
three cases.

{\bf Case 1.} $\dim[x_1,x_2,Ax_1,Ax_2]=4$.

Writing $Ax_1=x_3$ and $Ax_2=x_4$, by Hahn-Banach Theorem, there
exist $f_i\in X^*$ such that $f_i(x_j)=\delta_{ij}$ (the Kronecker's
symbol), $i,j=1,2,3,4$. Let $g_1=f_1+f_3$, $g_2=\alpha f_2+f_4$ with
$\alpha^{r+s-1}\neq1$, $|\alpha|=1$ and let $B=x_1\otimes g_1+x_2\otimes
g_2$; then rank$B=2$ and $\sigma_\pi(B)=\{1,\alpha\}$. It follows from
$B^r=x_1\otimes g_1+\alpha^{r-1}x_2\otimes g_2$ that $B^rAB^s=x_1\otimes
g_1+\alpha^{r+s-2}x_2\otimes g_2$. Then we get $\sigma_\pi(B^rAB^s)=\{1,\alpha^{r+s-1}\}$,
 a contradiction.

{\bf Case 2.} $\dim[x_1,x_2,Ax_1,Ax_2]=3$.

Since $Ax_1=x_3$ and $Ax_2=x_4$ are linearly independent, we have
$x_2=\lambda_1x_1+\lambda_2x_3+\lambda_3x_4$ for some scalars
$\lambda_1$, $\lambda_2$, $\lambda_3$. Note that $\lambda_3\not=0$
as $\dim[x_1,x_2,Ax_1,Ax_2]=3$. By Hahn-Banach Theorem, there exist
$f_i\in X^*$ such that $f_i(x_j)=\delta_{ij}$, $i,j=1,3,4$. Let $g
_1=f_1+f_3$, $g_2=\alpha f_4$ with $\alpha\neq0$ and let
$B=x_1\otimes g_1+x_2\otimes g_2$. Then $$B^r=x_1\otimes g_1+
(\lambda_1+\lambda_2)\sum_{i=2}^{r}(\alpha\lambda_3)^{i-2}x_1\otimes
g_2 +(\alpha\lambda_3)^{r-1}x_2\otimes g_2, \;\;\;r\geq2.$$ We
consider six subcases:

{\bf Subcase 1.} s=0.

For any $x=\gamma_1x_1+\gamma_2x_2$ with $\gamma_2\neq0$, we have $B^rAx=
x_1\otimes g_1Ax+(\lambda_1+\lambda_2)\sum_{i=2}^{r}(\alpha\lambda_3)^{i-2}x_1\otimes g_2Ax
+(\alpha\lambda_3)^{r-1}x_2\otimes g_2Ax$ and $B^rAx_1=x_1$.
Thus,  $B^rAx=\lambda x$ for some scalar $\lambda$ if and only if
$$\gamma_1+\gamma_2(\lambda_1+\lambda_2)\alpha\sum_{i=2}^{r}(\alpha\lambda_3)^{i-2}=\lambda\gamma_1\eqno(2.3)$$ and
$$\gamma_2\alpha^r\lambda_3^{r-1}=\lambda\gamma_2.\eqno(2.4)$$
By Eq.(2.4), $\lambda=\alpha^r\lambda_3^{r-1}\not=0$ as $\lambda_3\neq0$.
Take $\alpha$ such that $|\alpha|^{r}=|\lambda_3|^{1-r}$ and
$\alpha^r\lambda_3^{r-1}\not=1$. Then
$\lambda=\alpha^r\lambda_3^{r-1}$ and $|\lambda|=1$. Let $\beta=\sum_{i=2}^{r}(\alpha\lambda_3)^{i-2}\alpha$. Take $\gamma_1=1$ and
$\gamma_2=\frac{\lambda-1}{(\lambda_1+\lambda_2)\beta}$ if
$\lambda_1+\lambda_2\not=0$ and $\beta\neq0$; take $\gamma_1=0$ and $\gamma_2=1$ if
$\lambda_1+\lambda_2=0$ or $\beta=0.$ Then $\lambda$ satisfies both Eq.(2.3) and
Eq.(2.4) and hence $B^rA(\gamma_1x_1+\gamma_2x_2)=\lambda
(\gamma_1x_1+\gamma_2x_2)$. This implies that
$\sigma_\pi(B^rA)=\{1,\lambda\}$, a contradiction.

{\bf Subcase 2.} r=0.

For any $x=\gamma_1x_3+\gamma_2x_4$ with $\gamma_2\neq0$, we have $AB^sx=
x_3\otimes g_1x+(\lambda_1+\lambda_2)\sum_{i=2}^{s}(\alpha\lambda_3)^{i-2}x_3\otimes g_2x
+(\alpha\lambda_3)^{s-1}x_4\otimes g_2x$ and $AB^sx_3=x_3$.
Thus, $AB^sx=\lambda x$ for some scalar $\lambda$ if and only if
$$\gamma_1+\gamma_2(\lambda_1+\lambda_2)\alpha\sum_{i=2}^{s}(\alpha\lambda_3)^{i-2}=\lambda\gamma_1\eqno(2.5)$$ and
$$\gamma_2(\alpha\lambda_3)^{s-1}\alpha=\lambda\gamma_2.\eqno(2.6)$$
By Eq.(2.6), $\lambda=\alpha^s\lambda_3^{s-1}\not=0$ as $\lambda_3\neq0$.
Take $\alpha$ such that $|\alpha|^{s}=|\lambda_3|^{1-s}$ and
$\alpha^s\lambda_3^{s-1}\not=1$. Then
$\lambda=\alpha^s\lambda_3^{s-1}$ and $|\lambda|=1$. Let $\beta=\sum_{i=2}^{s}(\alpha\lambda_3)^{i-2}\alpha$. Take $\gamma_1=1$ and
$\gamma_2=\frac{\lambda-1}{(\lambda_1+\lambda_2)\beta}$ if
$\lambda_1+\lambda_2\not=0$ and $\beta\neq0$; take $\gamma_1=0$ and $\gamma_2=1$ if
$\lambda_1+\lambda_2=0$ or $\beta=0.$ Then $\lambda$ satisfies both Eq.(2.5) and
Eq.(2.6) and hence $AB^s(\gamma_1x_3+\gamma_2x_4)=\lambda
(\gamma_1x_3+\gamma_2x_4)$. This implies that
$\sigma_\pi(AB^s)=\{1,\lambda\}$, a contradiction.

{\bf Subcase 3.} r=s=1.

For any $x=\gamma_1x_1+\gamma_2x_2$ with $\gamma_2\neq0$, we have
$BABx=\gamma_1x_1+\gamma_2(\lambda_1+\lambda_2)x_1+\alpha^2\gamma_2\lambda_3x_2.$
Thus,  $BABx=\lambda x$ for some scalar $\lambda$ if and only if
$$\gamma_1+\gamma_2(\lambda_1+\lambda_2)=\lambda\gamma_1\eqno(2.7)$$ and
$$\alpha^2\gamma_2\lambda_3=\lambda\gamma_2.\eqno(2.8)$$
By Eq.(2.8), $\lambda=\alpha^2\lambda_3\not=0$ as $\lambda_3\neq0$.
Take $\alpha$ such that $|\alpha|^2=|\lambda_3|^{-1}$ and
$\alpha^2\lambda_3\not=1$. Then
$\lambda=\alpha^2\lambda_3$ and $|\lambda|=1$. Take $\gamma_1=1$ and
$\gamma_2=\frac{\lambda-1}{\lambda_1+\lambda_2}$ if
$\lambda_1+\lambda_2\not=0$; take $\gamma_1=0$ and $\gamma_2=1$ if
$\lambda_1+\lambda_2=0$. Then $\lambda$ satisfies both Eq.(2.7) and
Eq.(2.8) and hence $BAB(\gamma_1x_1+\gamma_2x_2)=\lambda
(\gamma_1x_1+\gamma_2x_2)$. This implies that
$\sigma_\pi(BAB)=\{1,\lambda\}$, a contradiction.

{\bf Subcase 4.} r=1.

For any $x=\gamma_1x_1+\gamma_2x_2$ with $\gamma_2\neq0$, we have $BAB^sx=
x_1\otimes g_1x+(\lambda_1+\lambda_2)\sum_{i=2}^{s}(\alpha\lambda_3)^{i-2}x_1\otimes g_2x
+(\alpha\lambda_3)^{s-1}\alpha x_2\otimes g_2x$ and $BAB^sx_1=x_1$.
Thus, $BAB^sx=\lambda x$ for some scalar $\lambda$ if and only if
$$\gamma_1+\gamma_2(\lambda_1+\lambda_2)(1+\sum_{i=2}^{s}(\alpha\lambda_3)^{i-1})=\lambda\gamma_1\eqno(2.9)$$ and
$$\gamma_2(\alpha\lambda_3)^{s}\alpha=\lambda\gamma_2.\eqno(2.10)$$
By Eq.(2.10), $\lambda=(\alpha\lambda_3)^{s}\alpha\not=0$ as $\lambda_3\neq0$.
Take $\alpha$ such that $|\alpha|^{s+1}=|\lambda_3|^{-s}$ and
$\alpha^{s+1}\lambda_3^s\not=1$. Then
$\lambda=\alpha^{s+1}\lambda_3^s$ and $|\lambda|=1$. Let $\beta=1+\sum_{i=2}^{s}(\alpha\lambda_3)^{i-1}$. Take $\gamma_1=1$ and
$\gamma_2=\frac{\lambda-1}{(\lambda_1+\lambda_2)\beta}$ if
$\lambda_1+\lambda_2\not=0$ and $\beta\neq0$; take $\gamma_1=0$ and $\gamma_2=1$ if
$\lambda_1+\lambda_2=0$ or $\beta=0.$ Then $\lambda$ satisfies both Eq.(2.9) and
Eq.(2.10) and hence $BAB^s(\gamma_1x_1+\gamma_2x_2)=\lambda
(\gamma_1x_1+\gamma_2x_2)$. This implies that
$\sigma_\pi(BAB^s)=\{1,\lambda\}$, a contradiction.

{\bf Subcase 5.} s=1.

For any $x=\gamma_1x_1+\gamma_2x_2$ with $\gamma_2\neq0$, we have
$$B^rABx= x_1\otimes
g_1x+(\lambda_1+\lambda_2)\alpha\sum_{i=2}^{r}(\alpha\lambda_3)^{i-2}x_1\otimes
g_2x +(\alpha\lambda_3)^{r-1}\alpha x_2\otimes g_2x$$ and
$B^rABx_1=x_1$. Thus, $B^rABx=\lambda x$ for some scalar $\lambda$
if and only if
$$\gamma_1+\gamma_2(\lambda_1+\lambda_2)(1+\alpha\sum_{i=2}^{r}(\alpha\lambda_3)^{i-1})=\lambda\gamma_1\eqno(2.11)$$ and
$$\gamma_2(\alpha\lambda_3)^{r}\alpha=\lambda\gamma_2.\eqno(2.12)$$
By Eq.(2.12), $\lambda=(\alpha\lambda_3)^{r}\alpha\not=0$ as $\lambda_3\neq0$.
Take $\alpha$ such that $|\alpha|^{r+1}=|\lambda_3|^{-r}$ and
$\alpha^{r+1}\lambda_3^r\not=1$. Then
$\lambda=\alpha^{r+1}\lambda_3^r$ and $|\lambda|=1$. Let $\beta=1+\sum_{i=2}^{r}(\alpha\lambda_3)^{i-1}\alpha$. Take $\gamma_1=1$ and
$\gamma_2=\frac{\lambda-1}{(\lambda_1+\lambda_2)\beta}$ if
$\lambda_1+\lambda_2\not=0$ and $\beta\neq0$; take $\gamma_1=0$ and $\gamma_2=1$ if
$\lambda_1+\lambda_2=0$ or $\beta=0.$ Then $\lambda$ satisfies both Eq.(2.11) and
Eq.(2.12) and hence $B^rAB(\gamma_1x_1+\gamma_2x_2)=\lambda
(\gamma_1x_1+\gamma_2x_2)$. This implies that
$\sigma_\pi(B^rAB)=\{1,\lambda\}$, a contradiction.

{\bf Subcase 6.} $s\geq2$, $r\geq2$.

For any $x=\gamma_1x_1+\gamma_2x_2$ with $\gamma_2\neq0$, we have
$$B^rAB^sx=\gamma_1x_1+\gamma_2(\lambda_1+\lambda_2)
(1+\sum_{j=2}^{s}(\alpha\lambda_3)^{j-1}+
\sum_{i=2}^{r}(\alpha\lambda_3)^{i+s-2}\alpha)x_1+\gamma_2(\alpha\lambda_3)^{r+s-1}\alpha x_2.$$
Thus,  $B^rAB^sx=\lambda x$ for some scalar $\lambda$ if and only if
$$\gamma_1+\gamma_2(\lambda_1+\lambda_2)(1+\sum_{j=2}^{s}(\alpha\lambda_3)^{j-1}+
\sum_{i=2}^{r}(\alpha\lambda_3)^{i+s-2}\alpha)=\lambda\gamma_1\eqno(2.13)$$ and
$$\gamma_2(\alpha\lambda_3)^{r+s-1}\alpha=\lambda\gamma_2.\eqno(2.14)$$
By Eq.(2.14), $\lambda=(\alpha\lambda_3)^{r+s-1}\alpha\not=0$ as
$\lambda_3\neq0$. Take $\alpha$ such that
$|\alpha|^{r+s}=|\lambda_3|^{1-r-s}$ and
$\alpha^{r+s}\lambda_3^{r+s-1}\not=1$. Then
$\lambda=\alpha^{r+s}\lambda_3^{r+s-1}$ and $|\lambda|=1$. Let
$\beta=1+\sum_{j=2}^{s}(\alpha\lambda_3)^{j-1}+
\sum_{i=2}^{r}(\alpha\lambda_3)^{i+s-2}\alpha$. Take $\gamma_1=1$
and $\gamma_2=\frac{\lambda-1}{(\lambda_1+\lambda_2)\beta}$ if
$\lambda_1+\lambda_2\not=0$ and $\beta\neq0$; take $\gamma_1=0$ and
$\gamma_2=1$ if $\lambda_1+\lambda_2=0$ or $\beta=0.$ Then $\lambda$
satisfies both Eq.(2.13) and Eq.(2.14) and hence
$B^rAB^s(\gamma_1x_1+\gamma_2x_2)=\lambda
(\gamma_1x_1+\gamma_2x_2)$. This implies that
$\sigma_\pi(B^rAB^s)=\{1,\lambda\}$, a contradiction.

 {\bf Case 3.} $\dim[x_1,x_2,Ax_1,Ax_2]=2$.

The above condition tells that in this case
$X_0=[x_1,x_2]=[Ax_1,Ax_2]$ is a $A$-invariant subspace of $X$. Let
$A_1$ be the restriction of $A$ to this subspace. It is invertible
and similar either to diag$(\alpha,\beta)$ with $\alpha\neq\beta$ or
to an upper triangular matrix with equal diagonal elements. In both
cases it is easy to construct $B_1$ such that
$\sigma_\pi(B_1^rA_1B_1^s)$ and hence, $\sigma_\pi(B^rAB^s)$ contains two
points.

The contradiction obtained in all cases imply that $A$ must have
rank one, as desired. \hfill$\Box$

Now let us give our proof of Theorem 2.2.

{\bf Proof of Theorem 2.2.} Since $B^rAB^s$ is a generalized
quasi-semi Jordan product of $A,B$, the ``if" part is true.

In the following we check the ``only if" part. Assume that the range
of $\Phi$ contains all operators of rank at most two and $\Phi$
satisfies Eq.(2.2). \if Since the case $r+s=1$ reduces to the case
of usual product, which is already discussed in \cite{TL}, we always
assume in the following that $r+s\geq 2$, and then Lemma 2.3 is
applicable.\fi

{\bf Claim 1.} For any $A\in\mathcal{A}_1$, $\Phi(A)=0$ if and only if $A=0$.

Let $\Phi(0)=B$. To prove $B=0$, assume, on the contrary, $B\neq0$;
then there exists a vector $x\in X_2$  such that $Bx\not=0$.

If $x$ and $Bx$ are linearly dependent, take $f\in X^*_2$ such that
$\langle Bx,f\rangle=1$, $\langle x,f\rangle\neq0$.
Let $T=x\otimes f$; then $T^rBT^s=(x\otimes f)^rB(x\otimes f)^s=(x\otimes f)^{r+s-1}$ and
$\sigma_\pi(T^rBT^s)=\{\langle x,f\rangle^{r+s-1}\}$.

If $x$ and $Bx$ are linearly independent, by Hahn-Banach Theorem,
there exist $f_1,f_2\in X^*_2$ such that $\langle
x,f_1\rangle=1$, $\langle x,f_2\rangle=0$, $\langle Bx,f_1\rangle=0$
and $\langle Bx,f_2\rangle=1$. Let $f=f_1+f_2$ and $T=x\otimes f$.
Then $\langle x,f\rangle=1$, $\langle Bx,f\rangle=1$ and
 $\sigma_\pi(T^rBT^s)=\sigma_\pi(T)=\{\langle x,f\rangle\}=\{1\}$.

Since the range of $\Phi$ contains all operators of rank at most
two, there exists $A\in\mathcal{A}_1$ such that $\Phi(A)=T$. Then
$$\{0\}=\sigma_\pi(A^r0A^s)=\sigma_\pi(\Phi(A)^r\Phi(0)\Phi(A)^s)=\sigma_\pi(T^rBT^s)=\{\langle
x,f\rangle^{r+s-1}\},$$ a contradiction. Hence we must have $B=0$.

Next we prove that  $\Phi(A)=0$ implies $A=0$. If $\Phi(A)=0$, then
we have
$$\sigma_\pi(B^rAB^s)=\sigma_\pi(\Phi(B)^r\Phi(A)\Phi(B)^s)=\{0\}$$
holds for all $B\in\mathcal{A}_1$, which forces that $A=0$.

 {\bf Claim 2.}
$\Phi$ preserves rank one operators in both directions.

Assume that rank$A=1$; then Claim 1 implies that $\Phi(A)\neq0$. For
any $B\in {\mathcal A}_1$, by Lemma 2.3,
$\sigma_\pi(\Phi(B)^r\Phi(A)\Phi(B)^s)=\sigma_\pi(B^rAB^s)$ is a singleton.
 Since the range of $\Phi$ contains all operators of rank at most
 two, for any $C\in{\mathcal A}_2$ with rank$(C)\leq 2$, $\sigma_\pi(C^r\Phi(A)C^s)$ is a singleton.
Applying Lemma 2.3 one has $\Phi(A)$ is of rank one. Conversely,
assume that $\Phi(A)$ is of rank one. Then, for any $B\in{\mathcal
A}_1$, Lemma 2.3 implies that
$\sigma_\pi(B^rAB^s)=\sigma_\pi(\Phi(B)^r\Phi(A)\Phi(B)^s)$ is a singleton.
Applying Lemma 2.3 again  one sees that $A$ is of rank one.

{\bf Claim 3.} $\Phi$ is linear and hence, by Claim 1, is injective.

We show first that $\Phi$ is additive. Note that, for any $A$ and
any rank-1 operator $x\otimes f$, we have
$$\sigma_\pi((x\otimes f)^rA(x\otimes f)^s)=\{\langle x,f\rangle\langle
Ax,f\rangle\}=\{{\rm Tr}((x\otimes f)^rA(x\otimes f)^s)\}.
\eqno(2.15)
$$

Let $A,B\in \mathcal{A}_1$ be arbitrary. For any $y\in X_2, g\in
X^*_2$ with $\langle y,g\rangle=1$, Claim 2 implies that there exist
$x\in X_1, f\in X^*_1$ such that $\Phi(x\otimes f)=y\otimes g$.
Then, by Eq.(2.15), we have
$$\begin{array}{rl} \{\langle\Phi(A+B)y, g\rangle\}
=&\sigma_\pi((y\otimes g)^r\Phi(A+B)(y\otimes g)^s)\\
=&\sigma_\pi((x\otimes f)^r(A+B)(x\otimes f)^s)\\
=&\{{\rm Tr}((x\otimes f)^r(A+B)(x\otimes f)^s)\}\\
=&\{{\rm Tr}((x\otimes f)^rA(x\otimes f)^s)+{\rm Tr}((x\otimes f)^rB(x\otimes f)^s)\}\\
=&\{{\rm Tr}((y\otimes g)^r \Phi(A) (y\otimes g)^s)+{\rm Tr}((y\otimes g)^r \Phi(B)(y\otimes g)^s)\}\\
=&\{{\rm Tr}((y\otimes g)^r(\Phi(A)+\Phi(B))(y\otimes g)^s)\}\\
=&\sigma_\pi((y\otimes g)^r(\Phi(A)+\Phi(B))(y\otimes g)^s)\\
=&\{\langle (\Phi(A)+\Phi(B))y, g\rangle\}.\end{array}$$  It follows
that
$$\langle\Phi(A+B)y,g\rangle=\langle(\Phi(A)+\Phi(B))y, g\rangle$$
holds for any $y\in X_2, g\in X^*_2$ with $\langle y,g\rangle=1$.
This entails $\Phi(A+B)=\Phi(A)+\Phi(B)$ and hence $\Phi$ is
additive. Similarly one can check that $\Phi$ is homogeneous. So
$\Phi$ is linear.

The claims 1-3 imply that  $\Phi$ is an injective linear map
preserving rank one operators in both directions.

Let us first consider the case  that $\dim X_1\geq 3$. Then, by
\cite{JH} the following claim is true.

{\bf Claim 4.} If $\dim X_1\geq 3$, then one of the following statements holds:

(i) There exist two linear bijections $T:X_1\rightarrow X_2$ and
$S:X^*_1\rightarrow X^*_2$ such that $\Phi(x\otimes f)=Tx\otimes Sf$ for
all rank one operators $x\otimes f\in \mathcal{A}_1$.

(ii) There exist two linear bijections $T:X^*_1\rightarrow X_2$ and
$S:X_1\rightarrow X^*_2$ such that $\Phi(x\otimes f)=Tf\otimes Sx$ for
all rank one operators $x\otimes f\in \mathcal{A}_1$.

{\bf Claim 5.} There exists a scalar $\lambda\in \mathbb{C}$ such
that  $\lambda^m=1$ with $m=r+s+1$ and, if (i) occurs in Claim $4$,
then $\langle Tx, Sf\rangle =\lambda \langle x,f\rangle$ holds for
all $x\in X_1$ and $f\in X^*_1$; if (ii) occurs in Claim $4$, then
$\langle Tf,Sx\rangle =\lambda \langle x,f\rangle$ holds for all
$x\in X_1$ and $f\in X^*_1$.

To check Claim 5, we first assume that the case (i) in Claim 4
occurs. Then, for any $x\in X_1$, $ f\in X^*_1$,
  we have $\sigma_\pi((x\otimes f)^m)=\{{\langle x, f\rangle}^m\}
  =\sigma_\pi((Tx\otimes Sf)^m)=\{{\langle Tx, Sf\rangle}^m\}$.
  So $\langle Tx, Sf\rangle=\lambda_{x,f}\langle x, f\rangle$ with $({\lambda_{x,f}})^m=1$.
   Especially, $\langle x, f\rangle=0\Leftrightarrow\langle Tx, Sf\rangle=0$.

Let $V_0=\{(x,f)\mid \langle x,f\rangle=0\}$, $V_t=\{(x,f)\mid\lambda_{x,f}
=e^{i\frac{2(t-1)\pi}{m}}\}$, $t=1,\ldots,m$.
 Then $\bigcup_{t=1}^m V_t=X_1\times X^*_1$ and $V_k\cap V_j=V_0$ if $k\not=j$,
$k,j=1,2,\ldots,m$. For $x_1,x_2\in X_1$ we denote by $[x_1,x_2]$ the linear
subspace spanned by $x_1$ and $x_2$.

 {\bf Assertion 1.} For any nonzero $x_1,x_2\in X_1, f\in X^*_1$,
 there exists some $k\in\{1,2,\ldots,m\}$ such that
$[x_1,x_2]\times[f]\subseteq V_k$.

We need only to show that we may take $\lambda_{x_1,f}$ and
$\lambda_{x_2,f}$ such that $\lambda_{x_1,f}=\lambda_{x_2,f}$. Consider the following
three cases.

{\bf Case 1.} $x_1,x_2$ are linearly dependent.

Assume that $x_2=\alpha x_1$; then $\alpha\not=0$ and
$\alpha\lambda_{x_1,f}\langle x_1, f\rangle=\alpha\langle Tx_1,
Sf\rangle=\langle Tx_2, Sf\rangle=\alpha\lambda_{x_2,f}\langle x_1,
f\rangle$. So we may take $\lambda_{x_1,f}$ and $\lambda _{x_2,f}$ such taht
$\lambda_{x_1,f}=\lambda _{x_2,f}$.

 {\bf Case 2.} $x_1,x_2$ are
linearly independent and at least one of $\langle x_i, f\rangle$,
$i=1,2$ is not zero.

In this case, for any $\alpha,\beta\in \mathbb{C}$ we have
$$\alpha\lambda_{\alpha,\beta}\langle x_1,
f\rangle+\beta\lambda_{\alpha,\beta}\langle x_2,f\rangle=\langle
T(\alpha x_1+\beta x_2), Sf\rangle=\alpha\lambda_{x_1,f}\langle
x_1,f\rangle+\beta\lambda_{x_2,f}\langle x_2,f\rangle,\eqno(2.16)$$
where $\lambda_{\alpha,\beta}=\lambda_{\alpha x_1+\beta x_2,f}.$ Let
$$\eta=\left(\begin{array}{c}
\lambda_{x_1,f}\langle x_1,f\rangle\\
\lambda_{x_2,f}\langle x_2,f\rangle
\end{array}\right),\
\eta_0=\left(\begin{array}{c}
\langle x_1,f\rangle\\
\langle x_2,f\rangle
\end{array}\right),\
\xi=\left(\begin{array}{c}
\alpha\\
\beta
\end{array}\right)\in \mathbb{C}^2.$$
Then Eq.(2.16) implies that
$$\langle\eta,\xi\rangle=\lambda_{\alpha,\beta}\langle\eta_0,\xi\rangle$$
holds for any $\xi\in\mathbb{C}^2$. It follows that
$\langle\eta,\xi\rangle=0\Leftrightarrow\langle\eta_0,\xi\rangle=0$.
So, as the vectors in ${\mathbb C}^2$, we must have $\eta=\gamma\eta_0$
for some scalar $\gamma$. Now it is clear that
$\lambda_{x_1,f}=\lambda_{x_2,f}$.

{\bf Case 3.}  $x_1,x_2$ are linearly independent and $\langle x_1,
f\rangle=\langle x_1,f\rangle=0$.

Then $\langle Tx_1,Sf\rangle=\lambda_{x_1,f}\langle x_1,f\rangle=0=\langle
Tx_2,Sf\rangle=\lambda_{x_2,f}\langle x_2,f\rangle$. In this case it
is clear that we can take $\lambda_{x_1,f}$ and $\lambda_{x_2,f}$
such that $\lambda_{x_1,f}=\lambda_{x_2,f}$.

Similar to the previous discussion, we have

{\bf Assertion 2.} For any nonzero $x\in X_1$, $f_1,f_2\in X^*_1$, there exists some
$k\in\{1,2,\ldots,m\}$ such that $[x]\times[f_1,f_2]\subseteq V_k$.

 {\bf Assertion 3.} There exists a scalar $\lambda\in\mathbb{C}$ with $\lambda^m=1$
such that $\lambda_{x,f}=\lambda$ for all $x\in X_1$ and $f\in X^*_1$.

For any $f_0\neq0$, there exist $x_0$ such that $\langle x_0,
f_0\rangle=1$. Then $\langle Tx_0, Sf_0\rangle=\lambda_{x_0,f_0}$
and $(x_0,f_0)\in V_{k_0}$ for some $k_0\in\{1,2,\ldots,m\}$.
So, by Assertion 1, for any $x\in X_1$, we have $[x,x_0]\times
[f_0]\subseteq V_{k_0}$, which implies that $X_1\times[f_0]\subseteq
V_{k_0}$. Similarly, by Assertion 2 one gets, for any $x_0\neq0$,
$[x_0]\times X^*_1\subseteq V_{k_0}$. Thus we obtain that $X_1\times
X^*_1=V_{k_0}$.

Hence, there exists a scalar $\lambda\in \mathbb{C}$ with $\lambda^m=1$ such
that $\lambda_{x,f}=\lambda$ for all $x$ and $f$, that is, $\langle
Tx, Sf\rangle =\lambda \langle x,f\rangle$ holds for all $x\in X_1$
and $f\in X^*_1$. So Assertion 3 is true.

This completes the proof  of Claim 5 for the case (i) of Claim 4.

If the case (ii) in Claim 4 occurs, by a similar argument one can
show that there exists a scalar $\lambda$ with $\lambda^m=1$ such
that $\langle Tf,Sx\rangle =\lambda \langle x,f\rangle$ holds for
all $x\in X_1$ and $f\in X^*_1$. Hence the last conclusion of Claim 5 is
also true.

{\bf Claim 6.} There exists a scalar $\lambda $ with $\lambda^m=1$
such that one of the followings holds:

(1) There exists an invertible operator $T\in {\mathcal B}(X_1, X_2)$
such that $\Phi(x\otimes f)=\lambda T(x\otimes f)T^{-1}$ for all
$x\otimes f\in \mathcal{A}_1$.

(2) $X_1$ and $X_2$ are reflexive, and there exists an invertible
operator $T\in {\mathcal B}(X^*_1, X_2)$ such that $\Phi(x\otimes
f)=\lambda T(x\otimes f)^* T^{-1}$ for all $x\otimes f\in
\mathcal{A}_1$.

Suppose that the case (i) of Claim 4 occurs. Then by Claim 5, there
exists a scalar $\lambda\in \mathbb{C}$ with $\lambda^m=1$ such that $\langle
Tx, Sf\rangle =\lambda \langle x,f\rangle$ holds for all $x\in X_1$
and $f\in X^*_1$. If $\{x_n\}\subset X_1$ is a sequence such that
$x_n\rightarrow x$ and $Tx_n\rightarrow y$ for some $x\in X_1$ and
$y\in X_2$ as $n\rightarrow\infty$, then, for any $f\in X^*_1$, we have
$$\langle y, Sf\rangle=\lim_{n\rightarrow\infty}\langle Tx_n,
Sf\rangle= \lim_{n\rightarrow\infty}\lambda\langle x_n,
f\rangle=\lambda\langle x,f\rangle=\langle Tx, Sf\rangle.$$ As $S$
is surjective we must have  $y=Tx$. So the bijection $T$ is a closed
operator and thus  a bounded invertible operator. Since $\langle Tx,
Sf\rangle=\langle x, T^*Sf\rangle=\lambda\langle x, f\rangle$ holds for all
$x\in X_1$ and $f\in X^*_1$, we
see that $T^*S=\lambda I$, that is $S=\lambda(T^*)^{-1}$. It follows
from the case (i) of Claim 4 that $\Phi(x\otimes f)=Tx\otimes Sf=\lambda
Tx\otimes (T^*)^{-1}f=\lambda T(x\otimes f)T^{-1}$ for any rank one
operator $x\otimes f$, i.e., the case (1) of Claim 6 holds.

Suppose that the case (ii) of Claim 4 occurs. Then by Claim 5, there
exists a scalar $\lambda\in \mathbb{C}$ with $\lambda^m=1$ such that $\langle
Tf, Sx\rangle =\lambda \langle x,f\rangle$ holds for all $x\in X_1$
and $f\in X^*_1$. Similar to the above argument we can check that both
$T$ and $S$ are bounded invertible operators with $S=\lambda
(T^*)^{-1}$. It follows that $\Phi(x\otimes f)=\lambda T(x\otimes
f)^*T^{-1}$ for any $x\otimes f$, obtaining that the case (2) of Claim
6 holds. Moreover, by \cite{JH}, in this case both $X_1$ and $X_2$ are
reflexive.

{\bf Claim 7.} The theorem is true for the case that $\dim X_1\geq 3$..

Assume that we have the case (1) of Claim 6.  Let $A\in \mathcal{A}_1$ be
arbitrary. For any $x\in X_1$ and $f\in X^*_1$ with $\langle x,
f\rangle=1$, we have
$$\begin{array}{rl}
\{\langle Ax,f\rangle\}=&\sigma_\pi((x\otimes f)^rA(x\otimes f)^s)\\
=&\sigma_\pi((\lambda T(x\otimes f)T^{-1})^r\Phi(A)(\lambda T(x\otimes f)T^{-1})^s)\\
=&\sigma_\pi(\frac{1}{\lambda}(x\otimes f)^rT^{-1}\Phi(A)T(x\otimes f)^s)\\
=&\{\langle\frac{1}{\lambda}T^{-1}\Phi(A)Tx,
f\rangle\}.\end{array}$$ This implies that $\Phi(A)=\lambda
TAT^{-1}$ for any $A\in \mathcal{A}_1$ and hence $\Phi$ has the form
(1) of Theorem 2.2.

A similar argument shows that if the case (2) of Claim 6 occurs then
$\Phi$ has the form given in (2) of Theorem 2.2.

This completes the proof of Theorem 2.2 for the case that $\dim
X_1\geq 3$.

In the rest of the proof we consider the case that $\dim X_1\leq 2$.
By the assumption on the range of $\Phi$, if $\dim X_1=1$, then $\dim
X_2=1$, and in this case there exists a scalar $\lambda$ with
$\lambda^m=1$ such that $\Phi (A)=\lambda A$ for all $A$. So the
theorem is true for the case that $\dim X_1=1$. Next we consider the
case that $\dim X_1=2$.

{\bf Claim 8.} The theorem is true for the case that $\dim X_1=2$.

By Claim 1 and Claim 3, $\Phi$ is a linear injection. By Claim 2,
$\Phi$  preserves rank one operators in both directions. As the
range of $\Phi$ contains all operators in ${\mathcal B}(X_2)$ of rank
at most two, we see that $\dim X_2=2$ since the range of $\Phi$ is a
4-dimensional subspace of ${\mathcal B}(X_2)$.
 So we can identify
$\mathcal{A}_1$ and $\mathcal{A}_2$ with $M_2=M_2({\mathbb C})$ as $\dim
X_1=\dim X_2=2$. For any rank one operator $x\otimes f\in{\mathcal
A}=M_2({\mathbb C})$, write $\Phi(x\otimes f)=y\otimes g$. Then
$$\{\langle x,f\rangle^m\}=\sigma _\pi ((x\otimes f)^m)=\sigma _\pi ((y\otimes
g)^m)=\{\langle y,g\rangle^m\}$$ and hence $\Phi(x\otimes f)$ is
nilpotent if and only if $x\otimes f$ is. By \cite[Corollary
2.5]{HC}, there exist a nonzero scalar $c\in\mathbb{C}$, a
nonsingular matrix $T\in M_2$ and a linear map
$\varphi:\mathbb{F}I\rightarrow M_2$ such that one of the following
statements holds:

(1) $\Phi(A)=cTAT^{-1}+\varphi(({\rm tr}(A)I)$ for all $A\in M_2$.

(2) $\Phi(A)=cTA^tT^{-1}+\varphi(({\rm tr}(A)I)$ for all $A\in M_2$.

We may assume that (1) holds. Otherwise, replace $\Phi$ by the map
$A\mapsto\Phi(A^t)$. We may further assume that $T=I$. If this is not the case,
replace $\Phi$ by the map $A\mapsto T^{-1}\Phi(A)T$. So, without
loss of the generality, we may assume that
$$\Phi(A)=cA+\varphi({\rm
tr}(A)I)=cA+{\rm tr}(A)\varphi(I)$$ for all $A\in M_2$.

 Write
$\varphi(I)=\left(\begin{array}{cc} t_{11} & t_{12}\\ t_{21} &
t_{22} \end{array}\right)$. As $\Phi$ preserves rank one matrices in
both directions, for any $x=\left(\begin{array}{c} \xi_1\\ \xi_2
\end{array}\right)$ and $f=(\eta_1 \ \eta_2)$, the determinant of
$$\Phi(x\otimes f)= c\left(\begin{array}{cc} \xi_1\eta_1 & \xi_1\eta_2\\ \xi_2\eta_1&
\xi_2\eta_2
\end{array}\right)+(\xi_1\eta_1+\xi_2\eta_2)\left(\begin{array}{cc}
t_{11} & t_{12}\\ t_{21} & t_{22} \end{array}\right)$$ is a zero
function in $(\xi_1,\xi_2,\eta_1,\eta_2)$. That is, we have
$$ \begin{array}{rl}  0\equiv &
(t_{11}t_{22}-t_{12}t_{21}+ct_{22})\xi_1^2\eta_1^2
+(t_{11}t_{22}-t_{12}t_{21}+ct_{11}) \xi_2^2\eta_2^2
\\
&+[c(t_{11}+t_{22})+2(t_{11}t_{22}-t_{12}t_{21})]\xi_1\xi_2\eta_1\eta_2
\\ & -ct_{12}(\xi_1\xi_2\eta_1^2+\xi_2^2\eta_1 \eta_2)-ct_{21}(\xi_1^2
\eta_1\eta_2+\xi_1 \xi_2 \eta_2^2)
\end{array}
$$
As $c\not=0$, it follows that $t_{12}=t_{21}=0$,
$t_{11}(t_{22}+c)=t_{22}(t_{11}+c)=0$. Thus, $t_{11}=t_{22}$ and
they take value $0$ or $-c$. Hence, $\varphi(I)=0$ or $-cI$.

If $\varphi(I)=0$, then $\Phi$ has the form (1).

If $\varphi(I)=-cI$, then $\Phi(I)=-cI$ and $(-c)^{r+s+1}=1$. It
follows that $\Phi(E_{11})=-cE_{22}$, $\Phi(E_{22})=-cE_{11}$,
$\Phi(E_{12})=cE_{12}$ and $\Phi(E_{21})=cE_{21}$. Since $\Phi$
preserves the peripheral spectrum of the generalized product, we
have $(-c)^{r+s+1}=1$. Then, for any $A=(a_{ij})\in M_2$ we have
$$\begin{array}{rl}\Phi(A)=&\left(\begin{array}{cc}
-ca_{22} & ca_{12}\\
ca_{21} & -ca_{11}
\end{array}\right)=-c\left(\begin{array}{cc}
a_{22} & -a_{12}\\
-a_{21} & a_{11}
\end{array}\right)\\=&-c\left(\begin{array}{cc}
0 & 1\\
-1 & 0
\end{array}\right)\left(\begin{array}{cc}
a_{11} & a_{12}\\
a_{21} & a_{22}
\end{array}\right)^t\left(\begin{array}{cc}
0 & 1\\
-1 & 0
\end{array}\right)^{-1}.\end{array}$$
Hence $\Phi$ has the form (2) of the theorem.

The proof of Theorem 2.2 is completed. \hfill$\Box$

{\bf Remark 2.4.} The assumption that the range of $\Phi$ contains
all operators of rank $\leq 2$ can not be omitted even for
case that $\dim X_1= 2$. To see this, assume that $\dim X_2=3$ and
consider the map $\Phi:M_2\rightarrow M_3$ defined by
$$\Phi(A)=\left(\begin{array}{cc}A&0\\0&0\end{array}\right)
$$
for every $A\in M_2$. It is clear that $\Phi$ preserves the
peripheral spectrum of the generalized products but $\Phi$ is
not the form stated in Theorem 2.1.

\section{The skew generalized products of operators on Hilbert spaces}

Let $\mathcal{A}_1$ and $\mathcal{A}_2$ be standard operator
algebras on complex Hilbert spaces $H$ and $K$, respectively. In
this section, we  characterize the maps from $\mathcal A_1$ into
$\mathcal A_2$ preserving the peripheral spectrum of the skew
generalized products.

Similar to the definition of generalized product, the skew
generalized product of operators on Hilbert space is defined as
follows. Fix a positive integer $k$ and a finite sequence
$(i_1,i_2,\ldots,i_m)$ such that
$\{i_1,i_2,\ldots,i_m\}=\{1,2,\ldots,k\}$ and there is an $i_p$ not
equal to $i_q$ for all other $q$. A skew generalized product for
operators $T_1,\ldots,T_k$ is defined by   $$T_1\diamond
T_2\diamond\cdots\diamond T_k=T_{i_1\cdots
}T_{i_{p-1}}T_{i_p}^*T_{i_{p+1}}\cdots T_{i_m}. \eqno(3.1)$$ The
definition of Eq.(3.1)   covers the usual skew product $T_1T_2^*$,
skew Jordan semi-triple $BA^*B$ and the skew triple one
$\{T_1,T_2,T_3\}=T_1T_2^*T_3$, etc.. Furthermore, if $(i_1,i_2,
\ldots, i_m)$ is symmetrical with respect to $i_p$, the above skew
generalized product is said to be semi Jordan. For instance,
$T_1\diamond T_2\diamond T_3=T_2T_3^2T_1^*T_3^2T_2$ is a skew
generalized semi Jordan product. Similar to Definition 1.1,
$T_1\diamond T_2\diamond\cdots\diamond T_k$ is called a skew
generalized quasi-semi Jordan product if Eq.(1.2) is true.

 For any unitary operator $U:H\rightarrow K$
with $U{\mathcal A}_1U^*\subseteq{\mathcal A_2}$, it is clear that
the map $A\mapsto UAU^*$  preserves the peripheral spectrum of any
skew generalized products. If $U{\mathcal
A}_1^tU^*\subseteq{\mathcal A_2}$,   the map $A\mapsto UA^tU^*$
preserves the peripheral spectrum of any skew generalized quasi-semi
Jordan products, here $A^t$ stands for the transpose of $A$ in an
arbitrary but fixed orthonormal basis of $H$. Also observe that, in
the case that $m$ is even, the map $\Phi$ preserves the peripheral
spectrum of the skew generalized products if and only if $-\Phi$
does.

The following result says that the converse is also true.

{\bf Theorem 3.1.} {\it Let ${\mathcal A}_1$ and ${\mathcal A}_2$ be
standard operator algebras on complex Hilbert spaces $H$ and $K$,
respectively. Consider the skew product $T_1\diamond\cdots\diamond
T_k$ defined in Eq.(3.1) with width $m$. Assume that
$\Phi:\mathcal{A}_1\rightarrow\mathcal{A}_2$ is a map the range of
which contains all operators of rank at most two. Then $\Phi$
satisfies
$$\sigma_\pi(\Phi(A_1)\diamond\Phi(A_2)\diamond\cdots\diamond \Phi(A_k))
=\sigma_\pi(A_1\diamond  A_2\diamond\cdots\diamond A_k) \eqno(3.2)$$
for all $A_1,A_2\ldots, A_k\in{\mathcal A}_1$  if and only if  there
exist a unitary operator $U\in {\mathcal B}(H, K)$ and a scalar
$c\in\{-1,1\}$ such that either}

(1) {\it $\Phi (A)=cUAU^*$ for every $A\in\mathcal{A}_1$; or }

(2) {\it $\Phi (A)=cUA^tU^*$ for every $A\in\mathcal{A}_1$ if  the
skew generalized product is of quasi-semi Jordan. \\  Here $A^t$ is
the transpose of $A$ with respect to an arbitrary but fixed
orthonormal basis of $H$. Moreover, $c=1$ whenever $m$ is odd.}

Theorem 3.1 clearly follows from the special case of $k=2$ below, by
considering $A_{i_p}=A$ and all other $A_{i_q}=B$.

{\bf Theorem 3.2.} {\it Let $\mathcal{A}_1$ and $\mathcal{A}_2$ be
standard operator algebras on complex Hilbert spaces $H$ and $K$,
respectively. Assume that
$\Phi:\mathcal{A}_1\rightarrow\mathcal{A}_2$ is a map the range of
which contains all operators of rank at most two, and $r,s$ are
nonnegative integers with $r+s\geq 1$. Then $\Phi$ satisfies
$$\sigma_\pi(B^rA^*B^s)=\sigma_\pi(\Phi(B)^r\Phi(A)^*\Phi(B)^s)\quad\mbox{\it for all }\
A,B\in\mathcal{A}_1 \eqno(3.3)$$ if and only if there exist a
unitary operator $U\in {\mathcal B}(H, K)$ and a scalar
$c\in\{-1,1\}$ such that $\Phi (A)=cUAU^*$ for every
$A\in\mathcal{A}_1$ or $\Phi (A)=cUA^tU^*$ for every
$A\in\mathcal{A}_1$. Moreover, $c=1$ for the case $r+s$ is even.
Here $A^t$ is the transpose of $A$ with respect to an arbitrary but
fixed orthonormal basis of $H$.}

Theorem 3.2 can be proved by a similar approach as Theorem 2.2 with
some necessary modifications. There is another simpler approach if
we assume that ${\mathcal A}_1$ is unital and we give its detail
blow.

As rank$(A)=1$ if and only if rank$(A^*)=1$, the following lemma is
immediate from Lemma 2.3.

{\bf Lemma 3.3.} {\it Let $\mathcal{A}$ be a standard operator
algebra  on a complex Hilbert space  $H$ and $r,s$ be nonnegative
integers with $r+s\geq1$. For a nonzero operator $A\in\mathcal{A}$,
the following statements are equivalent.}

(1) {\it $A$ is of rank one.}

(2) {\it For any $B\in\mathcal{A}, \sigma_\pi(B^rA^*B^s)$ is a singleton.}

(3) {\it For any $B\in\mathcal{A}$  with {\rm rank}$(B)\leq2$,
$\sigma_\pi(B^rA^*B^s)$ is a singleton.}

The next lemma can be found in \cite[Corollary 1]{SDZD}.

{\bf Lemma 3.4.} {\it Let $A,B\in{\mathcal B}(H)$ be  nonzero
operators and $n\geq 2$ be an   integer. Then $\langle Ax,
x\rangle=\langle Bx, x\rangle ^n $ holds for any unit vector $x\in
H$ if and only if there exists a complex number $c$ such that
$A=c^nI$ and $B=cI$.}

Now we are in a position to give our proof of Theorem 3.2.

{\bf Proof of Theorem 3.2.} We need only check the ``only if" part.
Assume that $\Phi$ satisfies Eq.(3.3) and $I\in{\mathcal A}_1$. \if
One may show that $\Phi$ has the form (1) or (2) in Theorem 3.1 by a
similar approach as that in the proof of Theorem 2.2. But here we
give another approach.\fi

 {\bf Claim 1.} $\Phi$ preserves rank one
operators in both directions.

This is obvious by Lemma 3.3 and the assumption that the range of
$\Phi$ contains all operators of rank at most two.

{\bf Claim 2.} $\Phi (I)=I$ or $-I$. $\Phi(I)=-I$ may occur only if
$r+s$ is odd.

For any unit vector $y\in K$, there exist $u, h\in H$ such that
$\Phi (u\otimes h)=y\otimes y$. It follows from $$\sigma _\pi
((u\otimes h)^r(u\otimes h)^*(u\otimes h)^s)=\sigma _\pi ((y\otimes
y)^{r+s+1})=\{1\}$$ that $\|h\|^2\|u\|^2\langle u, h\rangle^{r+s-1}
=1$, and hence $\langle u, h\rangle^{r+s-1}>0$. Since $$\{\langle u,
h\rangle ^{r+s}\}=\sigma _\pi ((u\otimes h)^{r+s})=\sigma _\pi
((y\otimes y)^r\Phi (I)^*(y\otimes y)^s)=\{\langle y, \Phi
(I)y\rangle \}$$ and
$$\{\langle h, u\rangle \}=\sigma _\pi
((u\otimes h)^*)=\sigma _\pi (\Phi (I)^r(y\otimes y)\Phi
(I)^s)=\{\langle \Phi (I)^{r+s}y, y\rangle \},$$ we see that
$$\langle  \Phi
(I)y,y\rangle=\langle \Phi (I)^{r+s}y, y\rangle ^{r+s}$$ holds for
all unite vector $y\in K$. Then by Lemma 3.4, there exists a scalar
$c$ such that $\Phi(I)=c^{r+s}I$ and $\Phi(I)^{r+s}=cI$. Thus we get
$c^{(r+s)^2-1}=1$ and $c =\langle h,u\rangle $ whenever
$\Phi(u\otimes h)$ is a projection. It follows from
$\bar{c}^{r+s-1}=\langle u, h\rangle^{r+s-1}>0$ that $c^{r+s-1}=1$.
Note that
$\{1\}=\sigma_\pi(I)=\sigma_\pi(\Phi(I)^r\Phi(I)^*\Phi(I)^s)=\{c^{r+s}\bar{c}\}$.
This implies that $c^2=c^{r+s-1}=1$ and hence $c=\pm 1$, i.e.,
$\Phi(I)=\pm I$.

It is clear that $c=1$, i.e., $\Phi(I)=I$, if $r+s$ is even. So, if
$c=-1$, then $r+s$ must be odd.

Note that, in the case $r+s$ is odd,   $\Phi$ satisfies Eq.(3.3) if
and only if $-\Phi$ satisfies Eq.(3.3). Therefore, in the case
$\Phi(I)=-I$, one may
 replace $\Phi$ by $-\Phi$, and  still assume that $\Phi(I)=I$. So,
 without loss of generality, we assume that $\Phi(I)=I$ in the rest
 of the proof.

{\bf Claim 3.} If $\Phi(I)=I$, then $\Phi$ preserves rank one
projections in both directions, and, there exists a unitary or
conjugate unitary operator $U:H \rightarrow K$ such that $\Phi
(x\otimes x)=Ux\otimes Ux$ for every unit vector $x\in H$.

 Assume that $\Phi(I)=I$. Accept the same symbols as that in the proof of Claim 2,
 we see that, if $\Phi(u\otimes h)=y\otimes y$ is a projection, then $\langle u, h\rangle =1$ and
$\|h\|\|u\|=1$, which implies that $h=\alpha u$ with $\alpha
>0$. Let $x=\sqrt{\alpha }u$, then $\|x\|=1$ and $\Phi (x\otimes
x)=y\otimes y$. That is, for any unit vector $y\in K$, there exists
a unit vector $x\in H$ such that $\Phi (x\otimes x)=y\otimes y$.
Conversely, since $\Phi (I)=I$, it is easily checked that $\Phi$
preserves rank one projections. Hence $\Phi$ preserves rank one
projections in both directions.

 It follows that there exists a bijective map $T:H\rightarrow K$ such
that
$$\Phi(x\otimes x)=Tx\otimes Tx $$
for all unit vectors $x\in H$ and $T(\lambda x)=\lambda Tx$ for any
$\lambda\in {\mathbb C}$, $x\in H$. Then, for any unit vectors
$x,y\in H$, we have
$$\begin{array}{rl} &\{|\langle x,y\rangle|^2\}=\sigma _\pi ((x\otimes x)^r(y\otimes y)^*(x\otimes
x)^s)\\ =&\sigma _\pi ((Tx\otimes Tx)^r(Ty\otimes Ty)^*(Tx\otimes
Tx)^s)=\{|\langle Tx , Ty\rangle |^2\}.\end{array}$$ Hence
$$|\langle Tx, Ty\rangle |=|\langle x,y\rangle| \eqno(3.4)$$ holds
for all $x,y\in H$.

The Wigner's theorem \cite{U} states that every bijective map $T$
between Hilbert spaces $H$, $K$ satisfying Eq.(3.4) must has the
form $ Tx=\phi(x)Ux$ for any $x\in H$, where $\phi $ is a generally
nonlinear functional on $H$ satisfying $|\phi(x)|\equiv 1$ and $U$
is a unitary or  a conjugate unitary (i.e., anti-unitary) operator.
Thus, by Wigner's theorem, there exists a unitary or conjugate
unitary operator $U:H \rightarrow K$ such that $\Phi (x\otimes
x)=Ux\otimes Ux$ for every unit vector $x\in H$.

Now assume that $U$ is unitary. Let $A\in \mathcal A_1$ be
arbitrary. For any unit vector $x\in H$, since
$$\begin{array}{rl}\{\langle x, Ax\rangle \}=&\sigma _\pi ((x\otimes
x)^rA^*(x\otimes x)^s)\\ =& \sigma _\pi ((Ux\otimes Ux)^r\Phi
(A)^*(Ux\otimes Ux)^s)=\{\langle Ux, \Phi(A)Ux\rangle
\},\end{array}$$we have
$$\langle Ax, x\rangle =\langle \Phi (A)Ux, Ux\rangle \qquad
\mbox{ for all unit vectors }x\in H.$$ Hence we get $\Phi (A)=UAU^*$
for every $A\in \mathcal A_1$.

Assume that $U$ is conjugate unitary. Take arbitrarily an orthonormal basis
$\{e_i\}_{i\in \Lambda }$ of $H$ and define $J$ by $J (\sum _{i\in
\Lambda }\xi _ie_i)=\sum _{i\in \Lambda }\bar{\xi _i}e_i$. Then $J :
H\rightarrow H$ is conjugate unitary and $JA^*J =
A^t$, where $A^t$ is the transpose of $A$ in the
orthonormal basis $\{e_i\}_{i\in \Lambda }$ of $H$. Let $V = JU$.
Then $V:H\rightarrow K$ is unitary. For any $A\in \mathcal A_1$, we
have
$$\begin{array}{rl}\{\langle x,
Ax\rangle \}=&\sigma _\pi ((x\otimes x)^rA^*(x\otimes x)^s)=\sigma
_\pi ((Ux\otimes Ux)^r\Phi (A)^*(Ux\otimes Ux)^s)\\ =&\{\langle Ux,
\Phi (A)Ux\rangle \}=\{\langle x,U^*\Phi (A)^*Ux\rangle \}=\{\langle
x,V^*\Phi (A)^tVx\rangle \}\end{array}$$ holds for any $x\in H$,
which forces that $V^*\Phi (A)^tV=A$. Therefore, in this case
$\Phi(A)=VA^tV^*$ for all $A$. This completes the proof.\hfill$\Box$

\end{document}